\documentclass{amsart}
\usepackage{amsmath}
\usepackage{amssymb}
\usepackage{hyperref}
\pdfoutput=1
\hypersetup{%
bookmarks=false,%
pdffitwindow=true,%
pdftitle={REDUCED DISTANCE BASED AT SINGULAR TIME IN THE RICCI FLOW},%
pdfauthor={Joerg Enders},%
pdfnewwindow=true,%
colorlinks=true,%
citecolor=black,%
filecolor=black,%
linkcolor=black,%
urlcolor=black,%
}

\def\del{\partial}

\def\pf{\non {\bf Proof. }}
\newtheorem{theorem}{Theorem}[section]
\newtheorem*{theorem*}{Theorem}
\newtheorem{lemma}[theorem]{Lemma}
\newtheorem{proposition}[theorem]{Proposition}
\newtheorem{corollary}[theorem]{Corollary}

\theoremstyle{definition}
\newtheorem{definition}[theorem]{Definition}

\newtheorem{example}[theorem]{Example}

\theoremstyle{remark}
\newtheorem{remark}[theorem]{Remark}

\def\non{\noindent}
\def\ref{\href}

\numberwithin{equation}{section}

\begin{document}

\title{REDUCED DISTANCE BASED AT SINGULAR TIME IN THE RICCI FLOW}

\author{Joerg Enders}
\address{Department of Mathematics, Michigan State University, East Lansing, MI 48824}
\email{endersjo@msu.edu}
\thanks{The author was partially supported by NSF grant DMS-0604759.}


\date{Received September 2007.}


\begin{abstract}
In this paper, we define a reduced distance function based at a
point at the singular time $T<\infty$ of a Ricci flow. We also show
the monotonicity of the corresponding reduced volume based at time
T, with equality iff the Ricci flow is a gradient shrinking soliton.
Our curvature bound assumption is more general than the type I
condition.
\end{abstract}
\maketitle

\section{Introduction}
\label{section1}

In \cite{Hamilton1982}, Hamilton introduced the evolution of Riemannian manifolds by the Ricci flow
$$\frac{\del g(t)}{\del t}=-2Ric_{g(t)}.$$
\non Solutions to this nonlinear equation generally develop
singularities in finite time. An important tool when studying the
formation of singularities in geometric evolution equation are
monotone quantities. Perelman \cite{PerelmanI} made significant
progress to the understanding of singularities in the Ricci flow by
finding such quantities. We will consider the 'reduced volume' and
extend its applicability to singular time.

\smallskip

Let $(M,g(t))$ be a complete oriented maximal solution to the Ricci flow for $t\in[0,T)$ where $T<\infty$.
We further assume a new curvature assumption, which we call 'type A' (see section \ref{typeA} for the definition)
and which includes the well-known type I condition. For any smooth
$(p,t_0)$ in space-time $M\times[0,T),$ Perelman defines the 'reduced distance' $l_{p,t_0}.$ It is a function on
$M\times[0,t_0]$ and satisfies the differential inequality
$$-\frac{\del}{\del \bar{t}}l_{p,t_0}(q,\bar{t})- \Delta_{g(\bar{t})} l_{p,t_0}(q,\bar{t}) +|\nabla^{g(\bar{t})}l_{p,t_0}(q,\bar{t})|_{g(\bar{t})}^2-R_{g(\bar{t})} +\frac{n}{2(t_0-\bar{t})}\ge
0.$$
This inequality implies the monotonicity (in $\bar{t}$) of the 'reduced volume'
$$\tilde{V}_{p,t_0}(\bar{t}):=\int_M
(4\pi(t_0-\bar{t}))^{-\frac
  n2}e^{-l_{p,t_0}(q,\bar{t})}\,dvol_{g(\bar{t})}(q)$$
along the Ricci flow. It follows from Perelman's work that
$\tilde{V}_{p,t_0}(\bar{t})$ is constant in $\bar{t}$ if and only if
$(M,g(t))$ is isometric to Euclidean space with the flat
(non-evolving) metric. In the proof, the equality case corresponds
to a gradient shrinking soliton (see Definition \ref{defsoliton}).
However, since $t_0$ is a regular time, gradient shrinking solitons
other than flat Euclidean space cannot arise. This motivates the
generalization considered in the paper, which we outline now.

For $p\in M$ and a sequence $t_i\nearrow T,$ we consider $l_{p,t_i}$ as above. In our main Theorem \ref{thm}, we show
the subsequential convergence
$$l_{p,t_i}\xrightarrow{C^0_{loc}} l_{p,T},$$ where the differential inequality survives to the limit function $l_{p,T},$
i.e. it satisfies
$$-\frac{\del}{\del\bar{t}}l_{p,T}(q,\bar{t})-\Delta
l_{p,T}(q,\bar{t})+|\nabla
l_{p,T}(q,\bar{t})|^2-R_{g(\bar{t})}(q)+\frac{n}{2(T-\bar{t})}\ge 0.$$

\non We conclude in Corollary \ref{vsingmon} that the corresponding
reduced volume based at singular time $(p,T)$ satisfies the
monotonicity formula with equality if and only if we have a gradient
shrinking soliton (see precise statement in section \ref{vsing}).

\smallskip

We would like to point out, that in very recent work \cite{Naber},
Naber obtained very similar results as we present in this paper.
While our type A assumption is weaker than the type I assumption
(and $\kappa-$noncollapsedness) used in his work, he additionally
applies the results to show that the rescaling limit of a type I
singularity is a shrinking soliton. This is the application we had
in mind and were currently working out the details of. As soon as we
learned about the posting of \cite{Naber}, we decided to type up
this paper, in which we independently present the results on the
reduced distance and reduced volume monotonicity based at singular
time, that we had presented in talks only so far. This explains the
incompleteness of this paper in terms of applications. We would like
to add that \cite{Naber} also classifies noncompact 4-dimensional
shrinking solitons with bounded nonnegative curvature operator. As
the above discussion shows, classification of gradient shrinking
solitons is crucial in understanding singularities, and there are
many related recent results, e.g.
\cite{NiWallachII},\cite{NiWallachI},\cite{Kotschwar}.

\section{Background on Ricci flow} \label{section2} In this section
we present background material on Ricci flow and the
reduced distance needed in the subsequent sections.

\subsection{The equation} \label{equation}
Let $(M^n,g(t)),\,t\in[0,T),\,0<T<\infty$ be a 1-parameter family of
complete, oriented n-dimensional Riemannian manifolds with bounded
(sectional) curvature solving the equation
\begin{equation}\label{rf}
\frac{\del g(t)}{\del t}=-2 Ric_{g(t)},
\end{equation}
where $Ric_{g(t)}$ denotes the Ricci curvature tensor of the metric
$g(t).$ We call such a family a \textbf{Ricci flow on $[0,T)$}
throughout this paper.

It follows from work by Hamilton \cite{Hamilton1982}, DeTurck
\cite{DeTurck}, Shi \cite{ShiI}\cite{ShiII} as well as Chen and Zhu
\cite{ChenZhuUniqueness} that for a given complete Riemannian
manifold with bounded curvature $(M,g_0),$ there exists a unique
solution to the quasilinear second order weakly parabolic equation
(\ref{rf}) with $g(0)=g_0$ on a time interval $[0,T)$ and the
solution is maximal if and only if
$$\lim_{t\nearrow T} \max_M |Rm_{g(t)}|_{g(t)}=\infty,$$
where $Rm_{g(t)}$ denotes the full curvature tensor of $g(t).$

\subsection{Gradient shrinking solitons}\label{solitons}
We will now discuss special solutions to the Ricci flow (\ref{rf}):
Gradient shrinking solitons play an important role in the equality
case of monotone quantities, and hence in the study of
singularities.

Motivated by Perelman's perspective in \cite{PerelmanI} and
\cite{PerelmanII}, we make the following
\begin{definition}\label{defsoliton}
A Riemannian manifold $(M,g)$ is called a \textbf{gradient shrinking
soliton} if there exists a \textbf{potential function} $f: M
\rightarrow \mathbb{R}$ such that
\begin{equation}\label{solitonI}
Ric_g+\nabla^{g}\nabla^{g}f-\frac{1}{2}g=0,
\end{equation}
where $\nabla^{g}\nabla^{g}f$ denotes the covariant Hessian of $f.$
\end{definition}
If $\nabla^{g}f$ determines a complete vector field (e.g. if it is
bounded, or in particular if $M$ is compact), then as in Theorem 4.1
in \cite{ChowII} we get a corresponding Ricci flow $g(t)$ on
$(-\infty,T)$ with $g(T-1)=g$, called \textbf{gradient shrinking
soliton in canonical form}, in the following way: The 1-parameter
family of diffeomorphisms $\phi_t$ of $M$ generated by integrating
the vector field $\frac{1}{T-t}\nabla^{g}f$ with $\phi_{T-1}=id$
gives us
\begin{equation}\label{solitonII}
g(t)=(T-t)\phi_t^* g.
\end{equation}
With the potential function $f(t)=\phi_t^*f$ the solution $g(t)$
satisfies
\begin{equation}\label{solitonIII}
Ric_{g(t)}+\nabla^{g(t)}\nabla^{g(t)}f(t)-\frac{1}{2(T-t)}g(t)=0
\end{equation}
on $(-\infty,T).$ Also, $f(t)$ satisfies
\begin{equation}\label{dfdt}
\frac{\del f}{\del t}=|\nabla^{g(t)}f(t)|_{g(t)}^2.
\end{equation}

\begin{remark}\label{trickysol}
Note that if we have a Ricci flow $g(t)$ on $[0,T),$ which together
with a 1-parameter family of functions $f(t)$ satisfies
($\ref{solitonIII}),$ this clearly defines a gradient shrinking
soliton according to Definition \ref{defsoliton} by considering the
equation at $t=T-1.$ However, it might not be a gradient shrinking
soliton in canonical form. If however $\nabla^{g(T-1)}f(T-1)$ is
complete, we can conclude by uniqueness of solutions to the Ricci
flow that the corresponding Ricci flow in canonical form equals
$g(t)$ on $[T-1,T).$
\end{remark}

Solitons are 'self-similar' solutions, since they evolve only by
scaling and diffeomorphism. Their existence is to be expected due to
the diffeomorphism invariance of the Ricci flow equation (\ref{rf}).
Solitons can be regarded as generalized fixed points, i.e. fixed
points of the volume normalized Ricci flow on the space of metrics
modulo the diffeomorphism group.

By changing the minus sign in (\ref{solitonI}) to "+" (or dropping
the metric term), one can analogously define \textbf{gradient
expanding} (or \textbf{steady}) \textbf{solitons}. Also, if the
vector field generating the diffeomorphisms $\phi_t$ is not a
gradient vector field, one gets a more general notion of solitons.
We will skip a more detailed discussion, since those examples are
not directly relevant in this paper.

\begin{example}[\textbf{Einstein solutions}]\label{Einstein} Let $(M^n,g_0)$ be an Einstein manifold, 
  i.e. $g_0$ satisfies $Ric_{g_0}=\frac{R_{g_0}}{n}g_0,$
  where $R_{g_0}$ denotes the (constant) scalar curvature of $g_0.$ Then
  \begin{equation}\label{einstein}
  g(t)=(1-\frac{2}{n}R_{g_0}\,t)\,g_0
  \end{equation}
  is a solution the the Ricci flow, where $g(t)$ is Einstein for each
  $t$.
  It only changes by scaling and can be seen to be a fixed
  point of the volume normalized Ricci flow.
  If $R_{g_0}>0$ the solution exists on the time interval
  $(-\infty,\frac{n}{2R_{g_0}}).$ With $T:=\frac{n}{2R_{g_0}}$
  we then get from (\ref{einstein}) that \
  \begin{equation}\label{reinstein}
  R_{g(t)}=\frac{n}{2(T-t)} \qquad \text{ and } \qquad Ric_{g(t)}=\frac{1}{2(T-t)}g(t).
  \end{equation}
  It follows from Remark \ref{trickysol} that any
  Einstein solution with positive scalar curvature can be
  regarded as a gradient shrinking soliton with $f\equiv 0,$ which
  is in canonical form.
%
\end{example}

   \begin{example}\label{gaussian} Consider the non-evolving Ricci flow $(\mathbb{R}^n,g(t)=g_{\mathbb{R}^n}).$
  If we let $f(x,t)=\frac{|x|^2}{4(T-t)},$ we have $$\nabla\nabla f(t)=\frac{1}{2(T-t)}g_{\mathbb{R}^n},$$
  which makes flat Euclidean space into a gradient shrinking soliton, called \textbf{Gaussian
  soliton}. It is in canonical form and will arise in section \ref{pv}.
  \end{example}


  \begin{example} We can construct gradient shrinking solitons
  as products of Einstein solutions $N^{n-k}$ of positive scalar curvature with flat Euclidean space $\mathbb{R}^k.$
  It is an interesting question to ask when gradient shrinking solitons are of this form \cite{PetersenWylie}.
  \end{example}

  We will now discuss important equations satisfied by gradient
  shrinking solitons: Let $(M,g(t),f(t))$ be a gradient shrinking
  soliton in canonical form on $(-\infty,T)$. Tracing equation
  (\ref{solitonIII}) gives
  \begin{equation}\label{trace}
  R_{g(t)}+ \Delta_{g(t)} f(t) -\frac{n}{2(T-t)}=0.
  \end{equation}
  We conclude from equations (\ref{dfdt}) and (\ref{trace}) that
  \begin{equation}\label{fdiffeq}
  -\frac{\del f}{\del t}- \Delta_{g(t)} f(t) +|\nabla^{g(t)}f(t)|_{g(t)}^2-R_{g(t)} +\frac{n}{2(T-t)}=0.
  \end{equation}
  In Perelman's point of view \cite{PerelmanI}, if we let
  \begin{equation}\label{defu}
  u(x,t):=(4\pi(T-t))^{-\frac{n}{2}}e^{-f(x,t)},
  \end{equation}
  where $x\in M$ and $t\in (-\infty,T),$ then a straight forward computation shows that (\ref{fdiffeq}) is
  equivalent to
  \begin{equation}\label{udiffeq}
  \square^* u(x,t)=0.
  \end{equation}
  Here $\square^*=-\frac{\del}{\del t}-\Delta_{g(t)}+R_{g(t)}$ denotes
  the formal adjoint of the heat operator $\square=\frac{\del}{\del t}-\Delta_{g(t)}$
  under the Ricci flow. This observation plays a key role
  in the proof of the equality case of the reduced volume monotonicity in Corollaries $\ref{vmon}$ and $\ref{vsingmon}$.

\smallskip

We will also need the following fact: If $(M,g,f)$ is a gradient shrinking soliton, then it follows from the contracted
second Bianchi identity and (\ref{trace}) that there exists a constant $C\in\mathbb{R}$ such that
\begin{equation}\label{gradconst}
R_g+|\nabla^{g}f|^2_g-f=C.
\end{equation}

\subsection{Perelman's reduced distance} \label{pl} In this
section we will briefly discuss Perelman's reduced distance for the
Ricci flow. Contrary to \cite{PerelmanI}, we will use forward time
notation rather than backward time, only since it will come more
natural when in sections \ref{section3} and \ref{vsing} we consider
a sequence of different base-times. The monotonicity of Perelman's
reduced volume will then be described in section \ref{pv}. While the
results are due to Perelman, there are by now several references
detailing his work, e.g. \cite{SesumTianWang2004},
\cite{KleinerLott}, \cite{YeI}\cite{YeII}, \cite{CaoZhu},
\cite{MorganTian}, \cite{ChowIII}.

\begin{definition}
Let $(M^n,g(t))$ be a Ricci flow on $[0,T).$ For any curve
$\gamma:[\bar{t},t_0]\rightarrow M,$ where $0<\bar{t}<t_0<T,$ we
define the \textbf{$\mathcal{L}-$length} of $\gamma$ by
$$\mathcal{L}(\gamma):=\int_{\bar{t}}^{t_0}\sqrt{t_0-t}\bigl(|\dot{\gamma}(t)|_{g(t)}^2+R_{g(t)}(\gamma(t))\bigr)dt.$$
\end{definition}
\non At a curve $\gamma$ the first variation of $\mathcal{L}$ with
fixed endpoints is given by
\begin{equation}\label{variation}
\delta_Y \mathcal{L}(\gamma)=\langle 2 \sqrt{t_0-t}
\dot{\gamma}(t),Y(t)\rangle\,
\big|_{\bar{t}}^{t_0}\qquad\qquad\qquad\qquad\qquad\qquad\qquad
\end{equation}
$$-\int_{\bar{t}}^{t_0} 2 \sqrt{t_0-t}\langle
\nabla^{g(t)}_{\dot{\gamma}(t)}\dot{\gamma}(t)-\frac{1}{2(t_0-t)}\dot{\gamma}(t)-2Ric_{g(t)}(\dot{\gamma}(t),\cdot)^\#-\frac12\nabla^{g(t)}R_{g(t)},Y(t)\rangle
dt.$$

\non $Y(t)$ is the variational vector field along $\gamma(t)$ with
$Y(\bar{t})=Y(t_0)=0$ and $^\#$ denotes the metric dual. The
Euler-Lagrange equation is called \textbf{$\mathcal{L}$-geodesic
equation} and given by
\begin{equation}\label{lgeI}
\nabla^{g(t)}_{\dot{\gamma}(t)}\dot{\gamma}(t)-\frac{1}{2(t_0-t)}\dot{\gamma}(t)-2Ric_{g(t)}(\dot{\gamma}(t),\cdot)^\#-\frac12\nabla^{g(t)}R_{g(t)}=0.
\end{equation}
Smooth solutions to (\ref{lgeI}) are called
\textbf{$\mathcal{L}$-geodesics.} Multiplying by $t_0-t,$ we can
rewrite (\ref{lgeI}) to get rid of the unbounded coefficient of the
second term:
\begin{equation}\label{lgeII}
\nabla^{g(t)}_{\sqrt{t_0-t}\dot{\gamma}(t)}(\sqrt{t_0-t}\,\dot{\gamma}(t))-2\sqrt{t_0-t}Ric_{g(t)}(\sqrt{t_0-t}\,\dot{\gamma}(t),\cdot)^\#-\frac12(t_0-t)\nabla^{g(t)}R_{g(t)}=0.
\end{equation}
Using the direct method in the calculus of variations after a change
of variables in $\mathcal{L}$, as well as the ODE (\ref{lgeII})
together with the curvature boundedness assumption and Shi's
derivative estimates \cite{ShiII}, one obtains the following

\begin{proposition}\label{lexists}
For any $(q,\bar{t})$ and $(p,t_0)$ with $p,q\in M$ and
$0<\bar{t}<t_0<T$ there exists a $\mathcal{L}-$minimizing
$\mathcal{L}-$goedesic $\gamma:[\bar{t},t_0]\rightarrow M$ with
$\gamma(\bar{t})=q$ and $\gamma(t_0)=p.$ Moreover, $\lim_{t\nearrow
t_0}\sqrt{t_0-t}\dot{\gamma}(t)$ exists.
\end{proposition}

\begin{definition}
For $(q,\bar{t})$ and $(p,t_0)$ as in the Proposition \ref{lexists}
we define
\begin{enumerate}
  \item \textbf{the $\mathcal{L}-$distance from $(q,\bar{t})$ to $(p,t_0)$} $$L_{p,t_0}(q,\bar{t}):=\inf\{\mathcal{L}(\gamma)\,|\,\gamma:[\bar{t},t_0]\rightarrow M,\gamma(\bar{t})=q,\gamma(t_0)=p\},$$
  \item \textbf{the reduced distance based at $(p,t_0)$}
  $$l_{p,t_0}(q,\bar{t}):=\frac{L_{p,t_0}(q,\bar{t})}{2\sqrt{t_0-\bar{t}}},$$
  \item and $$v_{p,t_0}(q,\bar{t}):=(4\pi(t_0-\bar{t}))^{-\frac
  n2}e^{-l_{p,t_0}(q,\bar{t})}.$$
\end{enumerate}
\end{definition}

\non We will fix $(p,t_0)\in M\times (0,T)$ and regard
$L_{p,t_0},\,l_{p,t_0}$ and $v_{p,t_0}$ as functions on space-time
$M\times (0,t_0).$

By studying the second variation of $\mathcal{L}$ and obtaining a
Laplacian comparison theorem for the reduced distance, Perelman
derives the following three differential (in)equalities:

\begin{theorem}\label{diffiethm}
Fix $(p,t_0)\in M\times (0,T).$ Then for all $(q,\bar{t})\in M\times
(0,t_0)$ we have
  \begin{equation}\label{diI}
  -\frac{\del}{\del \bar{t}}l_{p,t_0}(q,\bar{t})- \Delta_{g(\bar{t})} l_{p,t_0}(q,\bar{t}) +|\nabla^{g(\bar{t})}l_{p,t_0}(q,\bar{t})|_{g(\bar{t})}^2-R_{g(\bar{t})} +\frac{n}{2(t_0-\bar{t})}\ge
0,
 \end{equation}
 \begin{equation}\label{div}
  \text{ or equivalently }, \qquad\qquad\square^*v_{p,t_0}(q,\bar{t})\le 0.\qquad\qquad\qquad\qquad\qquad\qquad
  \end{equation}
  \begin{equation}\label{diII}
  -|\nabla^{g(\bar{t})}l_{p,t_0}(q,\bar{t})|_{g(\bar{t})}^2 + R_{g(\bar{t})}+\frac{l_{p,t_0}(q,\bar{t})-n}{t_0-\bar{t}}+ 2\Delta_{g(\bar{t})} l_{p,t_0}(q,\bar{t}) \le 0.
  \end{equation}
  \begin{equation}\label{diIII}
  -2\frac{\del}{\del \bar{t}}l_{p,t_0}(q,\bar{t})+|\nabla^{g(\bar{t})}l_{p,t_0}(q,\bar{t})|_{g(\bar{t})}^2 -R_{g(\bar{t})} +\frac{l_{p,t_0}(q,\bar{t})}{t_0-\bar{t}}= 0.
  \end{equation}
\end{theorem}

\begin{remark}\label{liprem}
It can be shown (see e.g. \cite{YeI}) that $L_{p,t_0}(q,\bar{t})$ is
locally Lipschitz on $M\times (0,t_0).$ Hence the same is true for
$l_{p,t_0}(q,\bar{t})$ and $v_{p,t_0}(q,\bar{t}),$ and the functions
are differentiable a.e. in $q$ and $\bar{t}.$ The (in)equalities in
Theorem \ref{diffiethm} therefore hold in the barrier sense, and in
particular in the sense of distributions. For (\ref{diI}) this means
that for any nonnegative $\phi\in C_{cpt}^\infty(M\times(0,t_0))$
\begin{eqnarray}\label{distdi}
&\int_0^T\int_M \Bigl(\bigl(-\frac{\del}{\del t}l_{p,t_0}
+|\nabla^{g(t)}l_{p,t_0}|_{g(t)}^2\hspace{-0.3cm}&-R_{g(t)}
+\frac{n}{2(t_0-t)}\bigr)\,\phi \\
\nonumber & &+
\nabla^{g(t)}l_{p,t_0}\cdot\nabla^{g(t)}\phi\Bigr)dvol_{g(t)} dt \ge
0.
\end{eqnarray}
\end{remark}

\subsection{Reduced volume and monotonicity}\label{pv}
We can now define Perelman's reduced volume.

\begin{definition}\label{rv}
Let $(M^n,g(t))$ be a Ricci flow on $[0,T),$ and $(p,t_0)\in
M\times(0,T).$ Then for each $\bar{t}\in (0,t_0)$ we define the
\textbf{reduced volume based at $(p,t_0)$}
$$\tilde{V}_{p,t_0}(\bar{t}):=\int_M
v_{p,t_0}(q,\bar{t})\,dvol_{g(\bar{t})}(q)=\int_M
\big(4\pi(t_0-\bar{t})\big)^{-\frac
  n2}e^{-l_{p,t_0}(q,\bar{t})}\,dvol_{g(\bar{t})}(q).$$
\end{definition}

\begin{example}
For the Gaussian soliton (example \ref{gaussian}) one computes
$$l_{p,t_0}(q,\bar{t})=\frac{|q-p|^2}{4(t_0-\bar{t})},$$
and hence
$$\tilde{V}_{p,t_0}(\bar{t})=\int_{\mathbb{R}^n}
\big(4\pi(t_0-\bar{t})\big)^{-\frac
  n2}e^{-\frac{|q-p|^2}{4(t_0-\bar{t})}}\,dq=1.$$
\end{example}

\smallskip

The monotonicity of the reduced volume along the Ricci flow is now
essentially a consequence of inequality (\ref{div}) in Theorem
\ref{diffiethm}:

\newpage

\begin{corollary}[Monotonicity of the reduced volume]\label{vmon}
Under the same assumptions as in Definition \ref{rv}, we have
\begin{enumerate}
  \item $\frac{d}{d\bar{t}}\tilde{V}_{p,t_0}(\bar{t})\ge 0,$
  \item $\lim_{\bar{t}\nearrow t_0}\tilde{V}_{p,t_0}(\bar{t})=1,$
  \item $\tilde{V}_{p,t_0}(\bar{t}_1)=\tilde{V}_{p,t_0}(\bar{t}_2)$
  for $0<\bar{t}_1<\bar{t}_2<t_0$ if and only if $(M,g(t))$ is isometric to the Gaussian soliton and
  $\tilde{V}_{p,t_0}(\bar{t})\equiv 1.$
\end{enumerate}
\end{corollary}

\section{Reduced distance based at singular time} \label{section3}

\subsection{Motivation} \label{motivation}
In the proof of Corollary \ref{vmon} (iii), one actually obtains
that if the reduced volume is constant, then $(M,g(t_0-1))$ is a
gradient shrinking soliton. Considering the corresponding canonical
form and using the fact that at time $t_0$ the curvature is bounded
(in fact zero), one concludes that the soliton must be the Gaussian
soliton. If for a maximal Ricci flow on $[0,T)$ we are able to base
the reduced distance and volume at singular time $(p,T)$, depending
on the base point $p\in M,$ we expect to get gradient shrinking
solitons other than the Gaussian soliton whenever this generalized
reduced volume is constant. We will prove this in Corollary
\ref{vsingmon} in section \ref{vsing}. Before that, in this section
we will define a reduced distance based a singular time.

\begin{remark}
Results of the type described above are known for other monotone
quantities in geometric evolution equations, e.g. for Perelman's
$\lambda$ and $\mu$ functionals for the Ricci flow \cite{PerelmanI}
or Huisken's monotonicity formula for the mean curvature flow
\cite{Huisken90}. Note also that the equality case of Harnack type
inequalities similarly identifies gradient expanding and steady
solitons, see e.g. \cite{Hamilton93}, \cite{ChenZhu00}, \cite{Ni02}.
\end{remark}

Let $(M,g(t))$ be a maximal Ricci flow on $[0,T).$ Let $t_i\nearrow
T$ and $p\in M.$ Then for all $(q,\bar{t})\in M\times (0,T)$ the
reduced distance $l_{p,t_i}(q,\bar{t})$ is defined for large enough
$i$ and the differential inequality (\ref{diI}) holds for each such
$i.$ This raises two questions:
\begin{enumerate}
\item Does there exist a good limit $l_{p,T}:=\lim_{t_i\nearrow T} l_{p,t_i}?$
\item Does the differential inequality (\ref{diI}) hold for $l_{p,T}$?
\end{enumerate}

\begin{example}
If $(M,g(t))$ is an Einstein solution on $[0,T)$ with $R(g(0))>0.$
For $p\in M$ and $t_i\nearrow T$ it follows from an explicit
computation (see e.g. \cite{ChowIII}) that
$$l_{p,t_i}(q,\bar{t})\rightarrow \frac{n}{2}=:l_{p,T}(q,\bar{t})$$
uniformly on $M\times [a,b]\subset M\times [0,T).$ Then (\ref{diI})
holds for the constant function $l_{p,T}$ because of equation
$(\ref{reinstein}$).
\end{example}

\begin{example}\label{lsingsoliton}
Let $(M,g(t),f(t))$ be a gradient shrinking soliton on $(-\infty,T)$
in canonical form and let $p\in M$ and $t_i\nearrow T.$ Then from
the discussion in section 7.3 in \cite{ChowIII} we can conclude that
for a subsequence
$$l_{p,t_i}(q,\bar{t})\rightarrow f(q,\bar{t}) + C=:l_{p,T}(q,\bar{t}),$$
where the convergence is uniform on compact subsets of
$M\times(0,T)$ and $C$ is the constant coming from equation
(\ref{gradconst}) satisfied by $(M,g(T-1),f(T-1))$. In particular,
the limit $l_{p,T}$ is independent of $p$ and $t_i.$
\end{example}

In general, to answer both questions above positively, we need to
mildly strengthen the bounded curvature assumption.

\subsection{Type A Ricci flows}
\label{typeA}

\begin{definition}
A Ricci flow $(M,g(t))$ on $[0,T)$ is said to be of \textbf{type A}
if there exist $C>0$ and $r\in[1,\frac{3}{2})$ such that for all
$t\in[0,T)$
$$|Rm_{g(t)}|_{g(t)}\le \frac{C}{(T-t)^r}.$$
\end{definition}

\begin{remark}
Note that for $r=1$ this is known as the \textbf{type I} condition.
Our type A assumption is weaker. From the maximum principle for
$|Rm_{g(t)}|_{g(t)}$ it follows that for a maximal Ricci flow on
$[0,T)$
$$max_M|Rm_{g(t)}|_{g(t)}\ge\frac{1}{8(T-t)},$$ which implies that
curvature blow-up with $r<1$ is impossible. On the other hand, the
type I condition is assumed to be generic. To our knowledge it is
not known whether there are maximal Ricci flows which are not of
type A. The only known example which is not of type I (i.e. type II)
is the degenerate neckpinch \cite{GuZhu}, but its curvature blow-up
rate is not known.
\end{remark}

\begin{example}\label{typeAex}
Let $(M,g,f)$ be a gradient shrinking soliton where $\nabla^g f$ is
a complete vector field. Then it follows from equation
(\ref{solitonII}) that the corresponding Ricci flow in canonical
form is of type A, in fact of type I.
\end{example}

\subsection{Main Theorem}

In this section, we state and prove the main

\begin{theorem}\label{thm}
Let $(M,g(t))$ be a Ricci flow on $[0,T)$ of type A. Also let $p\in
M$ and $t_i\nearrow T.$ Then there exists a locally Lipschitz
function
$$l_{p,T}:M\times (0,T)\rightarrow \mathbb{R},$$
which is a subsequential limit
$$l_{p,t_i}\xrightarrow{C^0_{loc}} l_{p,T}$$
and which satisfies the differential inequality analogous to
(\ref{diI}), i.e. for all $(q,\bar{t})\in M\times(0,T)$
\begin{equation}\label{dilpt}
-\frac{\del}{\del\bar{t}}l_{p,T}(q,\bar{t})-\Delta
l_{p,T}(q,\bar{t})+|\nabla
l_{p,T}(q,\bar{t})|^2-R_{g(\bar{t})}(q)+\frac{n}{2(T-\bar{t})}\ge 0
\end{equation}
holds in the sense of distributions. Equivalently,
\begin{equation}\label{divsing}
\square^* v_{p,T}\le 0,
\end{equation} where
$v_{p,T}(q,\bar{t}):=(4\pi(T-\bar{t}))^{-\frac{n}{2}}e^{-l_{p,T}(q,\bar{t})}.$
\end{theorem}

\non \textit{Remark.} This theorem has very recently been
independently obtained in \cite{Naber} under the type I assumption
(and $\kappa-$noncollapsedness) using very similar techniques.

\smallskip

Theorem \ref{thm} allows us to make the following
\begin{definition}
Under the assumptions of Theorem \ref{thm} we define
$$l_{p,T}:M\times (0,T)\rightarrow \mathbb{R}$$ to be a
\textbf{reduced distance for $(M,g(t))$ based at singular time
$(p,T).$}
\end{definition}
We now give the proof of Theorem \ref{thm}.

\smallskip

\pf To simplify notation let $l_i:=l_{p,t_i}(q,\bar{t})$ and
$L_i:=L_{p,t_i}(q,\bar{t}).$ The proof will be in 3 steps.

\smallskip

\non 1. First, we will derive a basic uniform bound on $l_i$ on
compact subsets $K=K_1\times [a,b]\subset M\times (0,T).$ By
definition of $l_i$ it suffices to show such a bound for
$L_i(q,\bar{t})$ on $K.$ Let $\eta:[0,1]\rightarrow M$ be a
$g(0)-$geodesic with $\eta(0)=q$ and $\eta(1)=p.$ Fix $k\in (b,T)$
and consider
$$\gamma(t):=\left\{\begin{array}{ll}
\eta(\frac{t-\bar{t}}{k-\bar{t}})& t\in[\bar{t},k]\\
p & t\in (k,t_i].
\end{array}\right.$$
Since $|\eta'(s)|_{g(0)}^2=c$ for a constant $c$, the uniform
equivalence of the metrics along the Ricci flow on $[0,k]$ yields a
constant $D$ such that $|\eta'(s)|_{g(t)}^2\le D.$ The type A
assumption implies that there exist constants $C$ and $r$ such that
$|R|\le\frac{C}{(T-t)^r}.$ Then we get the following estimate:
\begin{eqnarray}\label{Libound}
  |L_i(q,\bar{t})| &\le& \left|\int_{\bar{t}}^{t_i}\sqrt{t_i-t}\bigl(|\dot{\gamma}(t)|_{g(t)}^2+R_{g(t)}(\gamma(t))\bigr)dt\right| \\
   \nonumber &\le& \int_{\bar{t}}^{k}\frac{\sqrt{t_i-t}}{(k-\bar{t})^2}\Big|\eta'\Big(\frac{t-\bar{t}}{k-\bar{t}}\Big)\Big|_{g(t)}^2dt+C\int_{\bar{t}}^{T}\frac{\sqrt{t_i-t}}{(T-t)^r}dt \\
  \nonumber  &\le& \frac{D\sqrt{T}}{k-b}+ \frac{2C}{3-2r}T^{\frac32-r}=:E,
\end{eqnarray}
i.e. we have uniform bounds in $i$ on any given compact subset $K.$
\smallskip

\non 2. Next we derive uniform derivative bounds for $L_i$ on
compact subsets $K$ of space-time as above. We will first prove the
following

\begin{lemma}\label{lemma}
Under the assumptions of Theorem \ref{thm} let $\gamma_i(t)$ be an
$\mathcal{L}-$minimizing $\mathcal{L}-$geodesics from $(q,\bar{t})$
to $(p,t_i),$ where $(q,\bar{t})\in K.$ Then there exists a constant
$G$ independent of $i$ such that for all $t\in[\bar{t},k]$
$$|\sqrt{t_i-t}\,\dot{\gamma}_i(t)|^2_{g(t)}\le G.$$
\end{lemma}

\pf Denote by $V_i(t)=\sqrt{t_i-t}\,\dot{\gamma}_i(t).$ Using the
$\mathcal{L}-$geodesic equation (\ref{lgeII}) we compute

\begin{eqnarray}\label{odi}
 \nonumber\frac{d}{dt}|V_i(t)|_{g(t)}^2&=&-2Ric(V_i(t),V_i(t))+2\langle\nabla_{\dot{\gamma}_i(t)}V_i(t),V_i(t)\rangle_t\\
\nonumber &=&-2Ric(V_i(t),V_i(t))+\frac{2}{\sqrt{t_i-t}}\langle\nabla_{V_i(t)}V_i(t),V_i(t)\rangle_t\\
\nonumber &=&-2Ric(V_i(t),V_i(t))\\
\nonumber & &+\frac{2}{\sqrt{t_i-t}}\langle2\sqrt{t_i-t}Ric_{g(t)}(V_i(t),\cdot)^\#+\frac12(t_i-t)\nabla^{g(t)}R_{g(t)},V_i(t)\rangle_t\\
\nonumber &=&2Ric(V_i(t),V_i(t))+\sqrt{t_i-t}\,\langle\nabla^{g(t)}R_{g(t)},V_i(t)\rangle_t\\
&\le&2\frac{C_1}{(T-t)^r}|V_i(t)|^2_{g(t)}+\frac{C_2}{(T-t)^{r-\frac12}}|V_i(t)|_{g(t)},
\end{eqnarray}
where in the last inequality, since $t\ge\bar{t}>0,$ we used Shi's
derivative estimates \cite{ShiII} combined with the type A
assumption to bound $\nabla^{g(t)}R_{g(t)}.$ Note that $C_1, C_2$
are constants depending on the type A constant C, $n$ and $\bar{t},$
but are independent of $i.$

Before integrating this ordinary differential inequality to conclude
the proof of the lemma, we need to get uniform bounds on each
$V_i(t)$ for some $t$ in a compact set of time: By definition of
$\mathcal{L}$ and using again the type A assumption we estimate
\begin{eqnarray*}
  \int_{\bar{t}}^{t_i}\frac{1}{\sqrt{t_i-t}}|V_i(t)|^2_{g(t)}dt &=& \mathcal{L}(\gamma_i)-\int_{\bar{t}}^{t_i}\sqrt{t_i-t}\,R_{g(t)}(\gamma_i(t))dt\\
   &\le& \mathcal{L}(\gamma_i)+ \frac{2C}{3-2r}T^{\frac32-r}.
\end{eqnarray*}
Now the integral mean value theorem yields the existence of
$\hat{t}_i\in [\bar{t},k],$ such that

\begin{eqnarray*}
  \frac{1}{\sqrt{t_i-\hat{t}_i}}|V_i(\hat{t}_i)|^2_{g(\hat{t}_i)}&=& \frac{1}{k-\bar{t}}\int_{\bar{t}}^{k}\frac{1}{\sqrt{t_i-t}}|V_i(t)|^2_{g(t)}dt\\
  &\le&\frac{1}{k-\bar{t}}\int_{\bar{t}}^{t_i}\frac{1}{\sqrt{t_i-t}}|V_i(t)|^2_{g(t)}dt\\
   &\le&\frac{1}{k-\bar{t}}\big(\mathcal{L}(\gamma_i)+\frac{2C}{3-2r}T^{\frac32-r}\big),
\end{eqnarray*}
or equivalently
\begin{eqnarray*}
  |V_i(\hat{t}_i)|^2_{g(\hat{t}_i)}&=& \frac{\sqrt{t_i-\hat{t}_i}}{k-\bar{t}}\big(\mathcal{L}(\gamma_i)+\frac{2C}{3-2r}T^{\frac32-r}\big)\\
   &\le&\frac{\sqrt{T}}{k-b}\big(E+\frac{2C}{3-2r}T^{\frac32-r}\big)=:F,
\end{eqnarray*}
since by choice of $\gamma_i$ the bound (\ref{Libound}) holds for
$\mathcal{L}(\gamma_i).$

W.l.o.g. we can assume that $|V_i(t)|_{g(t)}\ge 1$ and estimate
(\ref{odi}) for $t\in[a,k]$ to get
$$\frac{d}{dt}|V_i(t)|_{g(t)}^2\le \left(2\frac{C_1}{(T-k)^r}+\frac{C_2}{(T-k)^{r-\frac12}}\right)|V_i(t)|^2_{g(t)}=C_3 |V_i(t)|^2_{g(t)}.$$
Integrating this, we conclude that for all $t\in[\bar{t},k]$
\begin{equation}\label{Vbounds}
|V_i(t)|^2_{g(t)}\le F e^{C_3(t-\hat{t}_i)}\le F e^{C_3T}=:G.
\end{equation}
This proves the Lemma.\qed

\smallskip

To get the gradient bounds for $L_i$ recall that from
(\ref{variation})
$$\nabla^{g(t)} L_i(q,\bar{t})=-2\sqrt{t_i-\bar{t}}\,\dot{\gamma}_i(\bar{t}),$$
so with Lemma \ref{lemma} we obtain for $(q,\bar{t})\in K$
\begin{equation}\label{gradbound}
    |\nabla^{g(t)} L_i(q,\bar{t})|_{g(\bar{t})}=2|\sqrt{t_i-\bar{t}}\,\dot{\gamma}_i(\bar{t})|_{g(\bar{t})}\le \sqrt{2G}.
\end{equation}

For the time derivative bound for $L_i$ we compute
\begin{eqnarray*}
  \frac{\del}{\del \bar{t}}
L_i(q,\bar{t}) &=& \frac{d}{d\bar{t}}L_i(q,\bar{t})-\langle
\nabla^{g(\bar{t})} L_i(q,\bar{t}),\dot{\gamma_i}(\bar{t})\rangle_{\bar{t}} \\
   &=&  - \sqrt{t_i-\bar{t}}\,\bigl(|\dot{\gamma_i}(\bar{t})|_{g(\bar{t})}^2+R_{g(\bar{t})}(\gamma(\bar{t}))\bigr)+2\sqrt{t_i-\bar{t}}\,|\dot{\gamma}(\bar{t})|_{g(\bar{t})}^2\\
   &=& \frac{1}{\sqrt{t_i-\bar{t}}}\,|\sqrt{t_i-\bar{t}}\,\dot{\gamma_i}(\bar{t})|_{g(\bar{t})}^2-\sqrt{t_i-\bar{t}}R_{g(\bar{t})}(\gamma(\bar{t}))
\end{eqnarray*}
Using the type A assumption and Lemma \ref{lemma} we get the time
derivative bound for $(q,\bar{t})\in K$
\begin{equation}\label{timebound}
    \big|\frac{\del}{\del \bar{t}}
L_i(q,\bar{t})\big|\le
\frac{G}{\sqrt{k-b}}+\frac{C}{(T-b)^{r-\frac12}}.
\end{equation}

\smallskip

Recall that by Remark \ref{liprem} each $l_i$ is locally Lipschitz
on $M\times (0,t_i)$. The above bounds show that $l_i$ are in fact
uniformly locally bounded and Lipschitz on $M\times(0,T).$ Hence
there exists a locally Lipschitz function
$$l_{p,T}:M\times (0,T)\rightarrow \mathbb{R}$$ such that a subsequence,
denoted $l_i$ again, converges to a $l_{p,T}$ in
$C^0_{loc}(M\times(0,T)).$ This proves the first part of the
theorem.

\smallskip

 \non 3. To prove that the differential inequality (\ref{dilpt}) holds
in the sense of distributions, we first note that $l_i\in
W^{1,2}_{loc}(M\times(0,T))$ and the bounds derived above imply
$$|l_i|_{W^{1,2}_{loc}(M\times(0,T))}<C.$$ W.l.o.g. we can assume that
$$l_i \rightharpoonup l_{p,T}$$ weakly in $W^{1,2}_{loc}(M\times(0,T)).$
This implies for all $(V,\psi)\in
W^{1,2}_{loc}(M\times(0,T),\mathbb{R}^{n+1})$

$$\int_0^T\int_M \nabla^{g(t)}l_i \cdot V +\frac{\del}{\del t}l_i\, \psi\, dvol_{g(t)}dt\rightarrow \int_0^T\int_M  \nabla^{g(t)}l_{p,T} \cdot V +\frac{\del}{\del t}l_{p,T}\, \psi\, dvol_{g(t)}dt$$

\non In particular, if for nonnegative $\phi\in
C^{\infty}_{cpt}(M\times(0,T)),$ we let $\psi=-\phi$ and
$V=\nabla^{g(t)}\phi,$ we get the distributional convergence

$$\int_0^T\int_M \nabla^{g(t)}l_i \cdot \nabla^{g(t)}\phi
-\frac{\del}{\del t}l_i\, \phi\, dvol_{g(t)}dt\qquad\qquad$$
$$\rightarrow
\int_0^T\int_M  \nabla^{g(t)}l_{p,T} \cdot \nabla^{g(t)}\phi
-\frac{\del}{\del t}l_{p,T}\, \phi\, dvol_{g(t)}dt.$$
\non Comparing with the distributional formulation $(\ref{distdi}),$
we see that we now only need to show that for all nonnegative
$\phi\in C^\infty_{cpt}(M\times(0,T))$
$$\int_0^T\int_M |\nabla^{g(t)}l_i|^2_{g(t)}\,\phi\, dvol_{g(t)}dt\rightarrow \int_0^T\int_M |\nabla^{g(t)}l_{p,T}|^2_{g(t)}\,\phi\, dvol_{g(t)}dt.$$
\non It suffices to show for each $t\in(0,T)$ and nonnegative $\phi\in C^\infty_{cpt}(M)$
$$\int_M |\nabla^{g(t)}l_i|^2_{g(t)}\,\phi\, dvol_{g(t)}\rightarrow \int_M |\nabla^{g(t)}l_{p,T}|^2_{g(t)}\,\phi\, dvol_{g(t)}.$$
\non Since this is weak $L^2$ convergence of $\sqrt{\phi}\,\nabla^{g(t)}l_i,$ we have
$$\int_M |\nabla^{g(t)}l_{p,T}|^2_{g(t)}\,\phi\, dvol_{g(t)}\le \liminf_{i\rightarrow\infty} \int_M |\nabla^{g(t)}l_i|^2_{g(t)}\,\phi\, dvol_{g(t)}.$$
\non We now show the other direction
$$\limsup_{i\rightarrow\infty} \int_M |\nabla^{g(t)}l_i|^2_{g(t)}\,\phi\, dvol_{g(t)}\le \int_M |\nabla^{g(t)}l_{p,T}|^2_{g(t)}\,\phi\, dvol_{g(t)}$$
using an argument similar to Lemma 9.21 in \cite{MorganTian}. We rewrite
\begin{eqnarray}\label{MTtrick}
\nonumber&&\limsup_{i\rightarrow\infty} \int_M \big(|\nabla^{g(t)}l_i|^2_{g(t)}-|\nabla^{g(t)}l_{p,T}|^2_{g(t)}\big)\,\phi\, dvol_{g(t)}\\
&=&\limsup_{i\rightarrow\infty}\bigl( \int_M \langle \nabla^{g(t)}(l_i-l_{p,T}),\phi\nabla^{g(t)}l_{p,T}\rangle_t\, dvol_{g(t)}\\
\nonumber&&\qquad+\int_M \langle
\nabla^{g(t)}(l_i-l_{p,T}),\phi\nabla^{g(t)}l_i\rangle_t\,
dvol_{g(t)}\bigr).
\end{eqnarray}
We can approximate $\phi\nabla^{g(t)}l_{p,T}$ by $V_j\in
C^\infty_{cpt}(M,\mathbb{R}^n)$ and conclude by weak $L^2$
convergence of $\nabla^{g(t)}l_i$ that the first integral goes to
zero as $i\rightarrow\infty$. For the second integral, we first use
the $C^0_{loc}-$ convergence of $l_i\rightarrow l_{p,T}$ to get a
sequence $\epsilon_i\searrow 0,$ such that on supp$(\phi)$ we have
$$l_{p,T}-l_i+\epsilon_i>0.$$ Then the second integral above equals
$$\limsup_{i\rightarrow\infty} \int_M \langle \nabla^{g(t)}(l_i-l_{p,T}-\epsilon_i),\phi\nabla^{g(t)}l_i\rangle_t\,dvol_{g(t)}.$$
We multiply Perelman's differential inequality (\ref{diII}) for $l_i$ by $\phi$ and write it in the distributional sense for a
nonnegative $\psi\in C^{\infty}_{cpt}(M):$
$$-\int_M\langle \nabla^{g(t)}(\psi\phi),\nabla^{g(t)}l_i\rangle_t\, dvol_{g(t)}\le \int_M\frac{\psi \phi}{2}\big(|\nabla^{g(t)}l_i|_{g(t)}^2 - R_{g(t)}-\frac{l_i-n}{t_i-t}\big)dvol_{g(t)}.$$
By approximation in $W^{1,2},$ we can take $\psi=l_{p,T}-l_i+\epsilon_i\ge 0$ to be only locally Lipschitz, i.e. conclude
$$\int_M\langle \nabla^{g(t)}\big((l_i-l_{p,T}-\epsilon_i)\phi\big),\nabla^{g(t)}l_i\rangle_t\, dvol_{g(t)}\qquad\qquad$$
$$\qquad\qquad\le \int_M\frac{(l_{p,T}-l_i+\epsilon_i) \phi}{2}\big(|\nabla^{g(t)}l_i|_{g(t)}^2 - R_{g(t)}-\frac{l_i-n}{t_i-t}\big)dvol_{g(t)}.$$
Since the right-hand integrand is bounded on supp$(\phi)$ and
$l_{p,T}-l_i+\epsilon_i\rightarrow 0$ uniformly, we obtain
$$\limsup_{i\rightarrow\infty}\int_M\langle \nabla^{g(t)}\big((l_i-l_{p,T}-\epsilon_i)\phi\big),\nabla^{g(t)}l_i\rangle_t\, dvol_{g(t)}\le 0.$$
Now inserting the characteristic function $\chi_{supp(\phi)},$ we
can rewrite
\begin{eqnarray*}
&&\int_M\langle
\nabla^{g(t)}\big((l_i-l_{p,T}-\epsilon_i)\phi\big),\nabla^{g(t)}l_i\rangle_t\,
dvol_{g(t)}\\
&=&\int_M\langle\nabla^{g(t)}(l_i-l_{p,T}-\epsilon_i)\phi,\nabla^{g(t)}l_i\rangle_t\,
dvol_{g(t)}\\
&&+\int_M\langle(l_i-l_{p,T}-\epsilon_i)\nabla^{g(t)}\phi,\chi_{supp(\phi)}\nabla^{g(t)}l_i\rangle_t\,
dvol_{g(t)}\\
&=&\int_M\langle\nabla^{g(t)}(l_i-l_{p,T}),\phi\nabla^{g(t)}l_i\rangle_t\,
dvol_{g(t)}\\
&&+\int_M\frac12(l_i-l_{p,T}-\epsilon_i)\big(|\nabla^{g(t)}\phi|_{g(t)}^2+|\chi_{supp(\phi)}\nabla^{g(t)}l_i|_{g(t)}^2\big)\,
dvol_{g(t)}.
\end{eqnarray*}
As before, the last integral tends to zero because of the uniform
convergence of $l_{p,T}-l_i+\epsilon_i\rightarrow 0.$ This implies
that the second term in (\ref{MTtrick}) satisfies
$$\limsup_{i\rightarrow\infty}\int_M\langle \nabla^{g(t)}(l_i-l_{p,T}),\phi\nabla^{g(t)}l_i\rangle_t\, dvol_{g(t)}\le 0$$
and finishes the proof. \qed


\section{Reduced volume and monotonicity based at singular time}
\label{vsing}

\begin{definition}\label{defvsing}
Let $(M,g(t))$ be a Ricci flow on $[0,T)$ of type A. Let $p\in M,$
$t_i\nearrow T,$ and $l_{p,T}$ and $v_{p,T}$ as in Theorem
\ref{thm}. Then we define a \textbf{reduced volume based a singular
time $(p,T)$} by
$$\tilde{V}_{p,T}(\bar{t}):=\int_M v_{p,T}(q,\bar{t}) dvol_{g(\bar{t})}(q)=\int_M\big((4\pi(T-\bar{t})\big)^{-\frac{n}{2}}e^{-l_{p,T}(q,\bar{t})}dvol_{g(\bar{t})}(q).$$
\end{definition}

\begin{remark}\label{Fatou}
The finiteness of $\tilde{V}_{p,T}(\bar{t})$ for any Ricci flow
$(M,g(t))$ and any fixed $\bar{t}\in(0,T)$ follows from Fatou's
lemma and the finiteness of Perelman's reduced volume as follows:
\begin{eqnarray*}
  \tilde{V}_{p,T}(\bar{t}) &=& \int_M\big((4\pi(T-\bar{t})\big)^{-\frac{n}{2}}e^{-\lim_{t_i\nearrow T} l_{p,t_i}(q,\bar{t})}dvol_{g(\bar{t})}(q) \\
   &=& \int_M\lim_{t_i\nearrow T}\left(\big((4\pi(T-\bar{t})\big)^{-\frac{n}{2}}e^{-l_{p,t_i}(q,\bar{t})}\right)dvol_{g(\bar{t})}(q) \\
   &\le& \liminf_{t_i\nearrow T} \underbrace{\tilde{V}_{p,t_i}(\bar{t})}_{\le 1}   \\
   &\le& 1.
\end{eqnarray*}
\end{remark}

Now we state the analogue of Perelman's monotonicity (Corollary
\ref{vmon}) for the reduced volume based at singular time.

\begin{corollary}[Monotonicity of the reduced volume based at singular time]\label{vsingmon} Under the
assumptions as in Definition \ref{defvsing} we have
\begin{enumerate}
  \item $\frac{d}{d\bar{t}}\tilde{V}_{p,T}(\bar{t})\ge 0,$
  \item $\lim_{\bar{t}\nearrow T}\tilde{V}_{p,T}(\bar{t})\le 1,$
  \item a) If $\tilde{V}_{p,T}(\bar{t}_1)=\tilde{V}_{p,T}(\bar{t}_2)$
  for $0<\bar{t}_1<\bar{t}_2<T,$ then $(M,g(T-1))$ is a gradient shrinking
soliton with potential function $l_{p,T}(\,\cdot\,,T-1).$

\non b) If $(M,g(t),f(t))$ is a gradient shrinking soliton in
canonical form and $p\in M,$ then any $\tilde{V}_{p,T}(\bar{t})$ is
constant in $\bar{t}.$
\end{enumerate}
\end{corollary}

\non \textit{Remark.} A similar statement has very recently also
been obtained by \cite{Naber} under the type I assumption (and
$\kappa-$noncollapsedness).

\smallskip

\pf (i) If $M$ is compact and $l_{p,T}$ is smooth the proof follows
directly from (\ref{divsing}) in Theorem \ref{thm}:
\begin{eqnarray}\label{prelim}
\nonumber\frac{d}{d\bar{t}}\tilde{V}_{p,T}(\bar{t})&=&\int_M
\big(\frac{\del}{\del\bar{t}}v_{p,T} - R \big)dvol_{g(\bar{t})}\\
&=&\int_M \big(\frac{\del}{\del\bar{t}}v_{p,T} +\Delta v_{p,T} - R\big)dvol_{g(\bar{t})}\\
\nonumber&=&\int_M -\square^*v_{p,T}dvol_{g(\bar{t})}\\
\nonumber&\ge& 0.
\end{eqnarray}

In the general case we need to justify the differentiation under the
integral and the adding in of the Laplacian term. For both
arguments, we need to bound
\begin{equation}\label{bound}\int_M
e^{-l_{p,T}}|l_{p,T}|dvol_{g(\bar{t})}<\infty\end{equation} for
fixed time $\bar{t}\in(0,T).$ From Remark \ref{Fatou} we know that
$$\int_M e^{-l_{p,T}}dvol_{g(\bar{t})}<\infty, \quad \text{ and
hence } \quad \int_M e^{-|l_{p,T}|}dvol_{g(\bar{t})}<\infty.$$ Since
the proof of the finiteness of Perelman's reduced volume $V_{p,t_i}$
relies on comparison of the reduced distance $l_{p,ti}$ with the
square of the distance function, we also get that
$$\int_M e^{-\frac12 l_{p,T}}dvol_{g(\bar{t})}<\infty, \quad \text{ and hence } \quad
\int_M e^{-\frac12|l_{p,T}|}dvol_{g(\bar{t})}<\infty.$$ Now let
$N:=\{q\in M\,|\,l_{p,T}(q,\bar{t})\ge 0\}.$ Then
$$\int_M
e^{-l_{p,T}}|l_{p,T}|dvol_{g(\bar{t})}=\int_{M-N}
e^{-l_{p,T}}|l_{p,T}|dvol_{g(\bar{t})}+\int_N
e^{-|l_{p,T}|}|l_{p,T}|dvol_{g(\bar{t})},$$ where the second term is
finite since $\frac12|l_{p,T}|\le e^{\frac12|l_{p,T}|}.$ The first
integral is bounded because the type A curvature bound yields an
upper bound on $|l_{p,T}|$ on $N$ by dropping the energy part in the
$\mathcal{L}-$functional. This proves (\ref{bound}).

The proof of Lemma \ref{lemma} implies that in fact there exist
constants $C_1$ and $C_2$ (depending on $\bar{t}$) such that
$$|\nabla^{g(\bar{t})}l_{p,T}|^2_{g(\bar{t})}\le C_1|l_{p,T}|_{g(\bar{t})}+C_2.$$
(Note that $L_{p,T}$ gets bounded on compact sets by the constant
$E$ in that proof.) This implies that
\begin{equation}\label{bound2}
\int_M e^{-l_{p,T}}|\nabla^{g(\bar{t})}l_{p,T}|_{g(\bar{t})}^2
dvol_{g(\bar{t})}<\infty.
\end{equation}

We lastly show that also
\begin{equation}\label{dtbound}
\int_M e^{-l_{p,T}}\big|\frac{\del}{\del\bar{t}}l_{p,T}\big|
dvol_{g(\bar{t})}<\infty.
\end{equation}
The differential equality analogous to (\ref{diIII})
$$-2\frac{\del}{\del
\bar{t}}l_{p,T}(q,\bar{t})+|\nabla^{g(\bar{t})}l_{p,T}(q,\bar{t})|_{g(\bar{t})}^2
-R_{g(\bar{t})} +\frac{l_{p,T}(q,\bar{t})}{T-\bar{t}}= 0$$ is
satisfied by $l_{p,T}$ for the same reasons that (\ref{dilpt}) in
Theorem \ref{thm} holds. We can use it to bound
$|\frac{\del}{\del\bar{t}}l_{p,T}|.$ Then $(\ref{dtbound})$ follows
from (\ref{bound}) and (\ref{bound2}) and the bounded curvature
assumption.

Now we can justify what corresponds to the adding in of the
Laplacian term in (\ref{prelim}) in the distributional setting. As
in $(\ref{distdi})$, Theorem \ref{thm} implies that for fixed
$\bar{t}\in(0,T)$ and for any nonnegative $\phi\in
C_{cpt}^\infty(M)$
\begin{equation*}
\int_M \Big(\big(-\frac{\del}{\del \bar{t}}l_{p,T}
+|\nabla^{g(\bar{t})}l_{p,T}|_{g(\bar{t})}^2-R_{g(\bar{t})}
+\frac{n}{2(T-\bar{t})}\big)\,\phi \\
+\nabla^{g(\bar{t})}l_{p,T}\cdot\nabla^{g(\bar{t})}\phi\Big)dvol_{g(\bar{t})}\ge
0.
\end{equation*}
Note that $$\nabla^{g(\bar{t})}(e^{-l_{p,T}})=-e^{-l_{p,T}}\nabla^{g(\bar{t})}l_{p,T},$$
so using the bounds for the integrals proved above and approximating
$\phi=e^{-l_{p,T}}$ by $\phi_j\in C_{cpt}^\infty(M),$ we conclude that
\begin{equation*}
\int_M \big(-\frac{\del}{\del \bar{t}}l_{p,T}
-R_{g(\bar{t})}+\frac{n}{2(T-\bar{t})}\big)e^{-l_{p,T}}dvol_{g(\bar{t})}\ge
0.
\end{equation*}
Up to the factor of $(4\pi(T-\bar{t}))^{-\frac{n}{2}}$ this is the
integral we obtain when differentiating the integrand in the reduced
volume based at $(p,T)$. To conclude (i), we only have to justify
the differentiation under the integral sign. We will do so using the
standard dominated convergence argument. We bound, at
$\bar{t}\in(0,T),$ the reduced volume integrand difference quotient
$D(q,\bar{t},h),$  rewritten in integral form, by an integrable
function:
$$D(q,\bar{t},h)=\frac{1}{h}\int_0^h \big(-\frac{\del}{\del \bar{t}}l_{p,T}(q,\bar{t}+s)
-R_{g(\bar{t}+s)}(q)\hspace{4.5cm}$$
$$\hspace{3cm}+\frac{n}{2(T-\bar{t}-s)}\big)(4\pi(T-\bar{t}-s))^{-\frac{n}{2}}e^{-l_{p,T}(q,\bar{t}+s)}\frac{dvol_{g(\bar{t}+s)}(q)}{dvol_{g(\bar{t})}(q)}ds$$
\non (compare e.g. \cite{ChowIII} for the argument in Perelman's
case) Now from Remark \ref{Fatou} and the bound (\ref{dtbound}) as
well as the exponential growth bound of the volume form under Ricci
flow, which makes the quotient
$$\frac{dvol_{g(\bar{t}+s)}(q)}{dvol_{g(\bar{t})}(q)}$$ bounded near
$\bar{t}$ (where we have bounded curvature), we conclude that
$D(q,\bar{t},h)$ is dominated by an integrable function for $h$ near
$0.$

\smallskip

\non(ii) This follows from (i) and Remark \ref{Fatou}.

\smallskip

\non(iii) a) From the initial computation in (i), which we
have justified for noncompact $M$ and in the distributional setting, we get
\begin{equation*}
0=\tilde{V}_{p,T}(\bar{t}_2)-\tilde{V}_{p,T}(\bar{t}_1)=\int_{\bar{t_1}}^{\bar{t_2}}\frac{d}{d\bar{t}}\tilde{V}_{p,T}(\bar{t})d\bar{t}=\int_{\bar{t_1}}^{\bar{t_2}}\int_M
-\square^*v_{p,T}dvol_{g(\bar{t})}d\bar{t}.
\end{equation*}
\non From (\ref{divsing}) we conclude that in the sense of distributions
\begin{equation}\label{v=0}
\square^*v_{p,T}\equiv 0,
\end{equation} and hence by parabolic regularity
$l_{p,T}$ is smooth. It follows from the proof of Theorem \ref{thm}
that (\ref{diIII}) holds for $l_{p,T}:$
\begin{equation}\label{diIIIsing}
-2\frac{\del}{\del\bar{t}}l_{p,T}(q,\bar{t})+|\nabla^{g(\bar{t})}l_{p,T}(q,\bar{t})|_{g(\bar{t})}^2
-R_{g(\bar{t})}+\frac{l_{p,T}(q,\bar{t})}{T-\bar{t}}=0.
\end{equation}
Combining (\ref{v=0}) and (\ref{diIIIsing}) we get
$$w_{p,T}:=\left((T-\bar{t})(2\Delta_{g(\bar{t})} l_{p,T}-|\nabla^{g(\bar{t})} l_{p,T}|_{g(\bar{t})}^2+R)+l_{p,T}-n\right)v_{p,T}\equiv0.$$
\non Hence
$$0\equiv\square^*w_{p,T}=-2(T-\bar{t})\left|Ric_{g(\bar{t})}+\nabla^{g(\bar{t})}\nabla^{g(\bar{t})}l_{p,T}-\frac{1}{2(T-\bar{t})}\right|^2_{g(\bar{t})}v_{p,T},$$
where the last equality follows from Perelman's Proposition 9.1
\cite{PerelmanI}. Because of $v_{p,T}>0$ we conclude that
$(M,g(T-1))$ is a gradient shrinking soliton with potential function
$l_{p,T}(\,\cdot\,,T-1).$

\smallskip

\non b) Let $p\in M.$ From Example \ref{typeAex} we get that the type A assumption is satisfied.
Then Example \ref{lsingsoliton} shows that there exists a constant $C$ such that
$$l_{p,T}(q,\bar{t})=f(q,\bar{t})+C,$$ and in particular $l_{p,T}$ is smooth. Hence the inequality
$(\ref{dilpt})$ becomes the equality (\ref{fdiffeq}), i.e.
$\square^*v_{p,T}=0.$ From the arguments in (i), we conclude
$$\frac{d}{d\bar{t}}\tilde{V}_{p,T}(\bar{t})=0$$ for all
$\bar{t}\in(0,T),$ which finishes the proof of Corollary
\ref{vsingmon}. \qed

\specialsection*{\textbf{Acknowledgements}}
The author would like to thank Professor Jon Wolfson for support, encouragement and many
discussions, as well as Professor Xiaodong Wang for many useful conversations.
The author also extends his gratitude to Professors Lei Ni and Bennett Chow for the
opportunity to visit UCSD and get helpful comments regarding this work.

\bibliographystyle{alpha}
\bibliography{c:/papers/mcf}

\end{document}